\documentclass[a4paper]{amsart}
\usepackage[margin=3.5cm]{geometry}
\usepackage[utf8]{inputenc}
\usepackage{amsmath,amsfonts,amssymb,amsthm, multirow}
\usepackage{comment, stackrel}
\usepackage{hyperref}
\usepackage{makecell}
\usepackage{hhline}
\usepackage{tikz}
\usepackage{tikz-cd}

\theoremstyle{plain}
\newtheorem{theorem}{Theorem}[section]

\newtheorem{proposition}[theorem]{Proposition}

\newtheorem{corollary}[theorem]{Corollary}
\newtheorem{lemma}[theorem]{Lemma}

\renewcommand{\phi}{\varphi}

\theoremstyle{definition}
\newtheorem{definition}[theorem]{Definition}

\newtheorem{example}[theorem]{Example}

\theoremstyle{remark}
\newtheorem{remark}[theorem]{Remark}

\newcommand{\C}{\mathbb{C}}                   
\renewcommand{\P}{\mathbb{P}}                 
\newcommand{\A}{\mathbb{A}}                   
\newcommand{\F}{\mathbb{F}}                   
\renewcommand{\O}{\mathcal{O}}

\newcommand{\CH}{\mathcal{H}}                 
\newcommand{\CT}{\mathcal{T}}                 
\newcommand{\CC}{\mathcal{C}}

\newcommand{\End}{\operatorname{End}}

\newcommand{\im}{\operatorname{im}}

\newcommand{\Spec}{\operatorname{Spec}}
\newcommand{\adj}{\operatorname{adj}}
\newcommand{\divv}{\operatorname{div}}
\newcommand{\ord}{\operatorname{ord}}
\newcommand{\rk}{\operatorname{rk}}

\newcommand{\Supp}{\operatorname{Supp}}
\newcommand{\Aut}{\operatorname{Aut}}
\newcommand{\Disc}{\operatorname{Disc}}

\newcommand{\br}[1]{\{#1\}}

\title{Ramification divisors and twisted Nambu brackets of Keller maps}

\author{Geoffrey Mboya\textsuperscript{1}}
\author{Stephen Mboya\textsuperscript{2,3}}
\author{Carringtone Kinyanjui\textsuperscript{4}}
\thanks{\textsuperscript{1} Maseno University, Kenya. Email: \href{mailto:mboya@maseno.ac.ke}{mboya@maseno.ac.ke}.}
\thanks{\textsuperscript{2}University of Nairobi, Kenya.}
\thanks{\textsuperscript{3}M\"alardalen University, Sweden.}
\thanks{\textsuperscript{4}\href{https://haki.africa/}{haki ai}.}

\date{Version 2}

\begin{document}
\begin{abstract}
For a polynomial map $\phi\colon\A^n\to\A^n$ we study the pushforward twist of the Jacobian Nambu bracket by $\phi^*$. The Hamiltonian fields of pulled back functions span the image of the adjugate of $d\phi$, so the twist detects the ramification divisor of $\phi$. For a Keller map $F$ the twist is removable exactly when $F$ is an automorphism; otherwise $F$ is not proper, and every projective extension of $F$ has a non-zero effective ramification divisor supported on the boundary, whose image is the non-properness set. For the explicit three-dimensional counterexample, whose finite cubic completion is described in~\cite{BGVai}, the twisted Hamiltonian fields move the marked root of a binary cubic and have polar divisor equal to the ramification divisor of the completion. Examples on weighted scrolls illustrate equivariance, degeneration to an automorphism, and an $A_1$ singularity of the fibre.
\end{abstract}
\maketitle

\section{Introduction}

\noindent The Jacobian bracket of $n$ polynomials in $n$ variables,
\[ \br{f_1,\ldots,f_n}=\det\left(\frac{\partial f_i}{\partial x_j}\right), \]
is the algebraic expression of the volume form $dx_1\wedge\cdots\wedge dx_n$; for $n=2$ it is the Poisson bracket $f_xg_y-f_yg_x$. It satisfies the Filippov identity and makes the polynomial ring an $n$-Lie algebra. A polynomial map $\phi\colon\A^n\to\A^n$ multiplies the bracket by its Jacobian determinant $J_\phi$ under pull-back. If $J_\phi$ is a non-zero constant, that is, if $\phi$ is a Keller map, then a suitable multiple $\sigma$ of $\phi^*$ is an endomorphism of the bracket, and the Yau twist $[f_1,\ldots,f_n]_\sigma=\sigma\br{f_1,\ldots,f_n}$ is a Hom-Nambu algebra with structure map $\sigma$~\cite{Yau,HLS,AMS}. The question of this paper is what this algebraic structure says about the geometry of $\phi$.

Two examples in dimension two show the two possible behaviours. For the shear $\phi(x,y)=(x+y^2,y)$ we have $J_\phi=1$ and the twisted bracket $[f,g]_\sigma=\phi^*\br{f,g}$ satisfies $[x,y]_\sigma=1$ and $[x,y^2]_\sigma=2y$. The shear is invertible, so $\sigma$ is invertible and $\br{f,g}=\sigma^{-1}[f,g]_\sigma$ can be recovered inside the algebra. We call such a twist removable. For the fold $\phi=(x^2,y)$ the Jacobian is $2x$ and the Hamiltonian fields of pulled back functions generate only the module $\langle\partial_x,x\partial_y\rangle$. On the ramification line $\{x=0\}$ this module drops to the kernel of $d\phi$, which is spanned by $\partial_x$.

The first main result, Theorem~\ref{thm:A}, says that this is the general picture. For any dominant polynomial map the Hamiltonian fields of pulled back functions generate the image of the adjugate of $d\phi$. Its cokernel is supported on the ramification divisor $R_\phi=\divv J_\phi$, with determinant divisor $(n-1)R_\phi$, and at points of corank one the fibre of the generated module is $\ker d\phi$. The result extends verbatim to relative maps over a base, Proposition~\ref{prop:relA}.

For a Keller map $F$ the ramification divisor is empty, and the twisted Hamiltonian fields are exactly the lifts of the Hamiltonian fields of the target, Proposition~\ref{prop:lift}. They form a flat connection on $F_*\O$, Proposition~\ref{prop:conn}. Locally the twist is invisible and its content is global. We prove that the twist is removable if and only if $F$ is an automorphism, Proposition~\ref{prop:untwist}, and that every Keller map that is not an automorphism is not proper, Proposition~\ref{prop:nonproper}. The cokernel $Q$ of $\C[F_1,\ldots,F_n]\subset\C[x_1,\ldots,x_n]$ is a flat bundle off the non-properness set whose monodromy is the permutation representation of the covering defined by $F$, Propositions~\ref{prop:Q} and~\ref{prop:mono}.

The non-properness set reappears at infinity. In Section~\ref{sec:infinity} we show that every projective extension $\bar F\colon\bar X\to\P^n$ of a Keller map that is not an automorphism has a non-zero effective ramification divisor $R^\infty_F$ supported on the boundary, Theorems~\ref{thm:C} and~\ref{thm:nonvanish}. The image of the boundary components that are not contracted to infinity is the non-properness set, Corollary~\ref{cor:Jel}; on the finite normalization of the target in the function field of the source~\cite{BGV} this is the statement that the non-properness set is the image of the boundary of the open embedding of $\A^n$ (Remark~\ref{rem:BGV}). For a family of Keller maps with a projective fibrewise compactification, the locus of fibres on which the twist is essential lies in the image of the support of a relative divisor, and in particular in a closed set, Theorem~\ref{thm:family}.

These statements are empty in dimension one and, conjecturally, in dimension two, where the Jacobian conjecture predicts that every Keller map is an automorphism. In dimension three and higher non-injective Keller maps exist; they are surveyed in~\cite{Gao}, and Proposition~\ref{prop:general} records what Sections~\ref{sec:twist}--\ref{sec:infinity} imply for each of them: the degree is at least two, the non-properness set is the locus where the fibre has fewer than $\deg F$ points, the monodromy is a transitive subgroup of the symmetric group, and the class of counterexamples is stable under composition with Keller automorphisms and under products with an identity. To make this concrete we study the explicit three-dimensional map announced by Alp\"oge in July 2026~\cite{Alpoge}, formally verified in~\cite{RHQ}, with $\det JF=-2$:
\begin{equation}\label{eq:Alpoge}
F_1=(1+xy)^3z+y^2(1+xy)(4+3xy),\quad
F_2=y+3x(1+xy)^2z+3xy^2(4+3xy),\quad
F_3=2x-3x^2y-x^3z.
\end{equation}
That $\A^3$ is the variety of binary cubic forms of the shape $u_3\mu^3-2\mu^2\eta+u_2\mu\eta^2-2u_1\eta^3$ with a marked simple root, that $F$ forgets the root, and the resulting stratification of the fibres by the discriminant surface $\Sigma$ and the triple-root curve $\Gamma$, are due to~\cite{Alpoge} and are worked out in the framework of~\cite{BGV} in~\cite{BGVai}; we recall them in Section~\ref{sec:cubic} in the form we use. From this, Theorems~\ref{thm:fibres} and~\ref{thm:twisted} read off the following. The cokernel $Q$ has rank $2$ off $\Sigma$ and is the bundle of the standard two-dimensional representation of $S_3$ there, and it vanishes over $\Sigma$. The twisted Hamiltonian fields are the vector fields moving the marked root, they have a simple pole along the ramification divisor $D$ of the finite cubic completion and no other pole, and the twisted ternary algebra contains an isomorphic copy of itself as a cubic extension. Proposition~\ref{prop:nu} shows that the coefficient of the boundary divisor over $D$ in the ramification divisor at infinity is exactly $1$ and that this divisor maps onto $\Sigma$.

Section~\ref{sec:families} treats families, and Section~\ref{sec:examples} gives three examples modelled on the weighted scrolls of~\cite{MS}. The first shows that the explicit counterexample is equivariant for a $\C^*$-action and lives on the scroll $\F(0,1,2,3\,|\,1^4)$. The second degenerates it along the weights of the scroll $\F(0,0,1,2\,|\,1^4)$ to an automorphism, with two of the three preimages escaping to a single point at infinity, so that the non-properness set is not stable under specialisation. The third is an $A_1$ degeneration of a double cover, where the absolute twisted Hamiltonian distribution does not see the node of the ramification curve while the relative one does.

\medskip
\noindent {\em Conventions.} We work over $\C$; varieties are integral, separated and of finite type. For $\phi\colon\A^n\to\A^n$ we write $\phi^*f=f\circ\phi$ and $d\phi=(\partial\phi_i/\partial x_j)$. A Keller map has $J_\phi=\det d\phi\in\C^\times$. A morphism of smooth varieties of equal dimension is a local isomorphism if its Jacobian determinant vanishes nowhere. By the Jacobian criterion~\cite[\S III.10]{Hartshorne}, \cite{Vakil} this is the same as being \'etale, and a local isomorphism is open with finite reduced fibres.

\section*{Acknowledgements}

\noindent We thank Adrian Vasiu for drawing our attention to~\cite{BGV} and~\cite{BGVai}, where the finite cubic completion of the explicit counterexample and the ramification divisor of the cubic cover are computed.

\section{Preliminaries}\label{sec:prelim}

\subsection{Nambu brackets}\label{sec:nambu}

Let $A=\C[x_1,\ldots,x_n]$ and $\omega=dx_1\wedge\cdots\wedge dx_n$. For $g_1,\ldots,g_{n-1}\in A$ the Hamiltonian field $X_g$ is defined by
\begin{equation}\label{eq:Xdef}
dg_1\wedge\cdots\wedge dg_{n-1}\wedge dh=X_g(h)\,\omega\qquad(h\in A).
\end{equation}
Equivalently $X_g(h)=\det(dg_1;\ldots;dg_{n-1};dh)$, with the differentials as rows. The Jacobian bracket is $\br{g_1,\ldots,g_{n-1},h}=X_g(h)$; it is skew-symmetric and a derivation in each argument~\cite{Takhtajan}. For example $\br{x,y,z}=1$ and $X_{x,y}=\partial_z$.

\begin{lemma}\label{lem:divfree}
Every Hamiltonian field is divergence free: $L_{X_g}\omega=0$.
\end{lemma}
\begin{proof}
Put $\eta=dg_1\wedge\cdots\wedge dg_{n-1}$ and $X=X_g$. Contracting the vanishing $(n+1)$-form $dh\wedge\omega$ with $X$ gives $X(h)\omega=dh\wedge\iota_X\omega$, while~\eqref{eq:Xdef} gives $X(h)\omega=(-1)^{n-1}dh\wedge\eta$. Hence the $(n-1)$-form $\beta=\iota_X\omega-(-1)^{n-1}\eta$ satisfies $dh\wedge\beta=0$ for every $h$. Taking $h=x_1,\ldots,x_n$ forces $\beta=0$, so $\iota_X\omega=\pm\eta$. This form is closed, being a product of exact forms, and Cartan's formula gives $L_X\omega=d\iota_X\omega=0$.
\end{proof}

\begin{proposition}[\cite{Filippov,Takhtajan}]\label{prop:Fil}
The Jacobian bracket satisfies the Filippov identity
\begin{equation}\label{eq:F}
\br{x_1,\ldots,x_{n-1},\br{y_1,\ldots,y_n}}=\sum_{i=1}^n\br{y_1,\ldots,\br{x_1,\ldots,x_{n-1},y_i},\ldots,y_n}.
\end{equation}
\end{proposition}
\begin{proof}
Let $\Lambda=\partial_1\wedge\cdots\wedge\partial_n$, so that $\br{h_1,\ldots,h_n}=\langle\Lambda,dh_1\wedge\cdots\wedge dh_n\rangle$. For a vector field $X$ the Leibniz rule for the Lie derivative gives
\[ X\br{h_1,\ldots,h_n}=\langle L_X\Lambda,dh_1\wedge\cdots\wedge dh_n\rangle+\sum_i\br{h_1,\ldots,Xh_i,\ldots,h_n}. \]
Since $\langle\Lambda,\omega\rangle=1$ we have $L_X\Lambda=-(\divv X)\Lambda$, and this vanishes for $X=X_{x_1,\ldots,x_{n-1}}$ by Lemma~\ref{lem:divfree}. Taking this $X$ and $h_i=y_i$ gives~\eqref{eq:F}.
\end{proof}

An algebra with a skew $n$-ary bracket satisfying~\eqref{eq:F} is an $n$-Lie algebra~\cite{Filippov}. An endomorphism of it is a linear map $\sigma$ with $\sigma\br{x_1,\ldots,x_n}=\br{\sigma x_1,\ldots,\sigma x_n}$.

\subsection{Hom-Nambu algebras and the Yau twist}\label{sec:yau}

\begin{definition}[\cite{AMS}]\label{def:Hom}
A multiplicative Hom-Nambu algebra is a triple $(L,[\cdot,\ldots,\cdot],\sigma)$ of a vector space, a skew $n$-linear bracket, and $\sigma\in\End(L)$ with $\sigma[x_1,\ldots,x_n]=[\sigma x_1,\ldots,\sigma x_n]$ and
\begin{equation}\label{eq:HF}
[\sigma x_1,\ldots,\sigma x_{n-1},[y_1,\ldots,y_n]]=\sum_{i=1}^n[\sigma y_1,\ldots,[x_1,\ldots,x_{n-1},y_i],\ldots,\sigma y_n].
\end{equation}
A morphism $(L,[\cdot],\sigma)\to(L',[\cdot]',\sigma')$ is a linear map $\psi$ with $\psi[x_1,\ldots,x_n]=[\psi x_1,\ldots,\psi x_n]'$ and $\psi\sigma=\sigma'\psi$.
\end{definition}

\begin{lemma}[Yau twist~\cite{Yau,HLS,AMS}]\label{lem:yau}
Let $(L,\br{\cdot})$ be an $n$-Lie algebra and $\sigma$ an endomorphism of it. Then $[x_1,\ldots,x_n]_\sigma:=\sigma\br{x_1,\ldots,x_n}$ makes $(L,[\cdot]_\sigma,\sigma)$ a multiplicative Hom-Nambu algebra.
\end{lemma}
\begin{proof}
Multiplicativity follows from $\sigma[x_1,\ldots,x_n]_\sigma=\sigma^2\br{x_1,\ldots,x_n}=\sigma\br{\sigma x_1,\ldots,\sigma x_n}=[\sigma x_1,\ldots,\sigma x_n]_\sigma$. Since $\sigma$ commutes with the bracket, the left side of~\eqref{eq:HF} equals $\sigma\br{\sigma x_1,\ldots,\sigma x_{n-1},\sigma\br{y_1,\ldots,y_n}}=\sigma^2\br{x_1,\ldots,x_{n-1},\br{y_1,\ldots,y_n}}$. The $i$-th term on the right equals $\sigma^2\br{y_1,\ldots,\br{x_1,\ldots,x_{n-1},y_i},\ldots,y_n}$. So~\eqref{eq:HF} is $\sigma^2$ applied to~\eqref{eq:F}.
\end{proof}

\subsection{Geometric facts}\label{sec:facts}

We use the following facts throughout.

\begin{enumerate}
\item[(G1)] A morphism of smooth varieties of equal dimension is a local isomorphism if and only if it is \'etale~\cite[\S III.10]{Hartshorne}, \cite{Vakil}. Such a morphism is open and quasi-finite with reduced fibres.
\item[(G2)] A proper quasi-finite morphism is finite~\cite[EGA\,III, 4.4.2]{EGA}. A finite morphism is proper in the classical topology, and properness is local on the target.
\item[(G3)] A proper local homeomorphism of connected locally compact Hausdorff spaces is a finite covering~\cite[Thm.\,4.22]{Forster}, and a connected covering of a simply connected space is a homeomorphism. The spaces $\C^n$ and $\P^n$ are simply connected.
\item[(G4)] An injective polynomial map $\A^n\to\A^n$ is an automorphism~\cite{Ax,BBR}.
\item[(G5)] Zariski's Main Theorem~\cite[EGA\,IV, 8.12.6]{EGA}: a quasi-finite separated morphism $X\to Y$ factors as an open immersion $X\hookrightarrow X^\nu$ followed by a finite morphism $X^\nu\to Y$. If $Y$ is normal and $X\to Y$ is dominant of degree $d$, then each fibre of $X^\nu\to Y$ has at most $d$ points~\cite[Ch.\,II, \S6.3]{Shafarevich}, hence so has each fibre of $X\to Y$. A quasi-finite birational morphism onto a normal variety is an open immersion. For a Keller map $F\colon\A^n\to\A^n$ one may take $X^\nu$ to be the normalization of the target in the function field of the source; this is the variety $X_F$ of~\cite{BGV}.
\item[(G6)] Riemann--Hurwitz~\cite[\S IV.2, Cor.\,2.4]{Hartshorne}: for a dominant, hence generically \'etale, morphism $g\colon X\to Y$ of smooth varieties of the same dimension $n$ the map $g^*\Omega^n_Y\to\Omega^n_X$ is injective. Its divisor of zeros is the ramification divisor $R_g\ge0$, and $K_X=g^*K_Y+R_g$.
\end{enumerate}

\section{The twist of a polynomial map}\label{sec:twist}

Throughout this section $A=\C[x_1,\ldots,x_n]$ carries the Jacobian bracket.

\begin{lemma}\label{lem:conf}
For every polynomial map $\phi\colon\A^n\to\A^n$ and all $f_1,\ldots,f_n\in A$,
\[ \br{\phi^*f_1,\ldots,\phi^*f_n}=J_\phi\cdot\phi^*\br{f_1,\ldots,f_n}. \]
\end{lemma}
\begin{proof}
By the chain rule the Jacobian matrix of $(f_i\circ\phi)_i$ is the product of the matrix $\big((\partial f_i/\partial u_k)\circ\phi\big)$ with $d\phi$. Take determinants.
\end{proof}

\begin{proposition}\label{prop:hom}
Let $F$ be a Keller map with $J_F=c$ and let $\kappa\in\C$ satisfy $\kappa^{n-1}c=1$. Then $\sigma=\kappa F^*$ is an endomorphism of $(A,\br{\cdot})$, and
\[ [f_1,\ldots,f_n]_\sigma:=\sigma\br{f_1,\ldots,f_n}=\kappa\,\br{f_1,\ldots,f_n}\circ F \]
makes $(A,[\cdot]_\sigma,\sigma)$ a multiplicative Hom-Nambu algebra.
\end{proposition}
\begin{proof}
By Lemma~\ref{lem:conf}, $\br{\sigma f_1,\ldots,\sigma f_n}=\kappa^nc\,F^*\br{f_1,\ldots,f_n}=\sigma\br{f_1,\ldots,f_n}$. Now apply Lemma~\ref{lem:yau} and Proposition~\ref{prop:Fil}.
\end{proof}

The normalisation is forced. For $\sigma=\lambda F^*$ the same computation gives $\br{\sigma f_1,\ldots,\sigma f_n}=\lambda^{n-1}c\,\sigma\br{f_1,\ldots,f_n}$, so $\sigma$ is an endomorphism of the bracket if and only if $\lambda^{n-1}c=1$.

\begin{example}\label{ex:shear}
For $n=2$ and $F(x,y)=(x+y^2,y)$ we have $c=\kappa=1$, $\sigma f(x,y)=f(x+y^2,y)$ and $[f,g]_\sigma=\sigma(f_xg_y-f_yg_x)$. Multiplicativity reads $[\sigma x,\sigma y]_\sigma=[x+y^2,y]_\sigma=\sigma(1)=1=\sigma[x,y]_\sigma$.
\end{example}

\begin{definition}\label{def:removable}
The twist of $(A,\br{\cdot})$ by $\sigma=\kappa F^*$ is removable if $\sigma\colon A\to A$ is bijective. Then $(A,[\cdot]_\sigma,\sigma)$ is the Yau twist by an automorphism and $\br{\cdot}=\sigma^{-1}[\cdot]_\sigma$ can be recovered inside $A$. Otherwise the twist is essential.
\end{definition}

\begin{proposition}\label{prop:untwist}
Let $F$ be a Keller map. Then $\sigma$ is injective, and the following are equivalent: (a) the twist is removable; (b) $F^*\colon A\to A$ is surjective; (c) $F$ is an automorphism. If $F$ is not an automorphism, then $\sigma(A)=\C[F_1,\ldots,F_n]$ is a proper subalgebra whose elements take the same value at any two points of a fibre of $F$.
\end{proposition}
\begin{proof}
The map $F$ is open and hence dominant by (G1), so $F^*$ is injective, and $\sigma(A)=F^*(A)$ because $\kappa\ne0$. This gives (a)$\Leftrightarrow$(b). If $F^*$ is surjective then $x_i=P_i(F_1,\ldots,F_n)$ for polynomials $P_i$, that is $P\circ F=\mathrm{id}$. Then $F$ is injective, hence an automorphism by (G4). Conversely an automorphism has a bijective $F^*$. Finally, $F(p)=F(p')$ implies $(f\circ F)(p)=(f\circ F)(p')$ for every $f$.
\end{proof}

\section{Ramification is degeneracy}\label{sec:ramification}

Let $\phi\colon\A^n\to\A^n$ be dominant with $n\ge2$. Identify vector fields on $\A^n$ with columns, so that $\CT_{\A^n}=\O^n$, and let $\CH_\phi\subset\CT_{\A^n}$ be the submodule generated by the fields $X_{\phi^*f_1,\ldots,\phi^*f_{n-1}}$ with $f_i\in A$. For a square matrix $M$ we have $M\adj M=\adj M\,M=(\det M)I$.

\begin{lemma}\label{lem:matrix}
Let $M\in M_n(R)$ over a commutative ring $R$, let $r_1,\ldots,r_{n-1}$ be rows, and let $v\in R^n$ be the column with $\det(r_1;\ldots;r_{n-1};t)=t\,v$ for all rows $t$. Then $\det(r_1M;\ldots;r_{n-1}M;s)=s\cdot\adj(M)\,v$ for every row $s$.
\end{lemma}
\begin{proof}
Both sides are polynomials with integer coefficients in the entries of $M$, $r_i$ and $s$, so it is enough to prove the identity over $\C$ with $\det M\neq0$. The left side is then $\det M\cdot\det(r_1;\ldots;r_{n-1};sM^{-1})=\det M\cdot sM^{-1}v=s\adj(M)v$.
\end{proof}

\begin{theorem}\label{thm:A}
Let $\phi\colon\A^n\to\A^n$ be dominant.
\begin{enumerate}
\item[(i)] $X_{\phi^*f_1,\ldots,\phi^*f_{n-1}}=\adj(d\phi)\cdot\big(X_{f_1,\ldots,f_{n-1}}\circ\phi\big)$, and $\CH_\phi=\adj(d\phi)(\O^n)$.
\item[(ii)] $\det\adj(d\phi)=J_\phi^{\,n-1}$. The cokernel of $\CH_\phi\hookrightarrow\CT_{\A^n}$ is supported on $R_\phi=\divv J_\phi$ and has determinant divisor $(n-1)R_\phi$.
\item[(iii)] Let $E_x\subset T_x\A^n$ be the span of the values at $x$ of the sections of $\CH_\phi$. Then $\dim E_x=n,1,0$ according as $\rk d\phi_x=n$, $n-1$, or $\le n-2$, and $E_x=\ker d\phi_x$ if $\rk d\phi_x=n-1$.
\end{enumerate}
\end{theorem}
\begin{proof}
(i) Put $M=d\phi$ and $g_i=f_i\circ\phi$. The differential $dg_i$ is the row $r_iM$ with $r_i=(df_i)\circ\phi$, and the column $v$ of Lemma~\ref{lem:matrix} is $X_f\circ\phi$. By~\eqref{eq:Xdef} and Lemma~\ref{lem:matrix} with $s=dh$ we get $X_g(h)=dh\cdot\adj(M)v$, which is the first claim. Take for $f$ the coordinate functions $(u_1,\ldots,\widehat{u_k},\ldots,u_n)$ of the target. Then $X_f=\pm\partial_{u_k}$, so $X_g$ is, up to sign, the $k$-th column of $\adj(d\phi)$. Hence $\CH_\phi$ contains $\adj(d\phi)(\O^n)$, and the formula for $X_g$ shows that it is contained in it.

(ii) The identity $\det\adj M=(\det M)^{n-1}$ is standard. Since $\phi$ is dominant it is generically \'etale in characteristic zero, so $J_\phi\ne0$ and $\adj(d\phi)$ is injective. It is invertible exactly where $J_\phi\ne0$. The determinant divisor of the cokernel of an injective endomorphism of $\O^n$ is the divisor of its determinant, here $(n-1)\divv J_\phi$.

(iii) The sections of $\CH_\phi$ are the columns of $\adj(d\phi)$ and their $\O$-combinations, so $E_x=\im\adj(d\phi_x)$. The adjugate of an $n\times n$ matrix of rank $n$, $n-1$, $\le n-2$ has rank $n$, $1$, $0$ respectively. If $\rk M=n-1$ then $M\adj M=0$ gives $\im\adj M\subseteq\ker M$, and both are lines.
\end{proof}

\begin{example}\label{ex:fold}
For $n=2$ and $\phi=(x^2,y)$ we have $d\phi=\mathrm{diag}(2x,1)$ and $\adj d\phi=\mathrm{diag}(1,2x)$. Hence $R_\phi=\{x=0\}$ and $\CH_\phi=\langle\partial_x,x\partial_y\rangle$, with cokernel $\O/(x)\,\partial_y$ of determinant divisor $R_\phi$, and $E=\langle\partial_x\rangle=\ker d\phi$ at $x=0$. The module $\CH_\phi$ is not involutive, because $[\partial_x,x\partial_y]=\partial_y\notin\CH_\phi$. So $\CH_\phi$ is a distribution only where $d\phi$ is invertible.
\end{example}

Theorem~\ref{thm:A} is the Riemann--Hurwitz formula in disguise. On $\A^n$ we have $\phi^*(du_1\wedge\cdots\wedge du_n)=J_\phi\,\omega$ with both forms nowhere zero, so $R_\phi=\divv J_\phi$ is the divisor of (G6). For $n=1$ and $\phi=x^e$ it is $(e-1)[0]$. For $n\ge2$ the same divisor, multiplied by $n-1$, is the degeneracy divisor of the twisted Nambu structure.

\section{Keller maps}\label{sec:keller}

Let $F\colon X\to Y$ be a local isomorphism of smooth affine varieties of dimension $n$, for instance a Keller map $\A^n\to\A^n$ by (G1). For a vector field $V$ on $Y$ there is a unique vector field $F^!V=(dF)^{-1}(V\circ F)$ on $X$ with $dF(F^!V)=V\circ F$. For a Keller map $(dF)^{-1}=c^{-1}\adj(dF)$, so $F^!V$ is polynomial whenever $V$ is.

\begin{lemma}\label{lem:lift}
(i) $(F^!V)(F^*g)=F^*(Vg)$. (ii) $F^!(gV)=F^*g\,F^!V$ and $[F^!V,F^!W]=F^![V,W]$.
\end{lemma}
\begin{proof}
(i) $(F^!V)(g\circ F)=dg_F\cdot dF\cdot F^!V=dg_F\cdot(V\circ F)=(Vg)\circ F$. (ii) The first identity is clear. Fields that are $F$-related have $F$-related bracket, so $dF[F^!V,F^!W]=[V,W]\circ F$, and a field is determined by its $F$-image because $dF$ is invertible.
\end{proof}

\subsection{Twisted Hamiltonian fields are lifts}

\begin{proposition}\label{prop:lift}
Let $F$ be a Keller map with $J_F=c$ and $\sigma=\kappa F^*$, where $\kappa^{n-1}c=1$. Then
\[ X_{F^*f_1,\ldots,F^*f_{n-1}}=c\,F^!X_{f_1,\ldots,f_{n-1}},\qquad X_{\sigma f_1,\ldots,\sigma f_{n-1}}=F^!X_{f_1,\ldots,f_{n-1}}, \]
and the fields $X_{\sigma f_1,\ldots,\sigma f_{n-1}}$ generate $\CT_{\A^n}$.
\end{proposition}
\begin{proof}
By Theorem~\ref{thm:A}(i), $X_{F^*f}=\adj(dF)(X_f\circ F)=c\,(dF)^{-1}(X_f\circ F)=c\,F^!X_f$. By multilinearity $X_{\sigma f}=\kappa^{n-1}X_{F^*f}=F^!X_f$. Since $\adj(dF)$ is invertible, $\CH_F=\O^n$.
\end{proof}

So for a Keller map the normalisation $\kappa^{n-1}c=1$ is exactly the one for which the twisted Hamiltonian fields are lifts without a residual scalar. The ramification divisor is empty and the twisted Hamiltonian distribution is the whole tangent bundle. Locally the twist is invisible, and its content is global.

\subsection{A connection and the cokernel}

Let $B=F^*\O(Y)\subset A=\O(X)$. For $h\in A$ put $\nabla_Vh=(F^!V)(h)$. By Lemma~\ref{lem:lift}, $\nabla_V(F^*g\cdot h)=F^*(Vg)h+F^*g\,\nabla_Vh$ and $[\nabla_V,\nabla_W]=\nabla_{[V,W]}$, so $\nabla$ is the canonical flat connection on $F_*\O_X$. In coordinates the fields $\xi_i:=F^!\partial_{u_i}$ satisfy $\xi_iF_k=\delta_{ik}$ and $[\xi_i,\xi_k]=0$.

\begin{proposition}\label{prop:conn}
(i) $A$ is a module over the Weyl algebra $\C\langle u_1,\ldots,u_n,\partial_1,\ldots,\partial_n\rangle$ via $u_k\mapsto F_k\cdot$ and $\partial_i\mapsto\xi_i$, and $B$ is the submodule generated by $1$. (ii) For a Keller map,
\[ [\sigma f_1,\ldots,\sigma f_{n-1},h]_\sigma=\sigma\big(\nabla_{X_{f_1,\ldots,f_{n-1}}}h\big). \]
\end{proposition}
\begin{proof}
(i) The Weyl relations are $[\xi_i,F_k]=\xi_iF_k=\delta_{ik}$ and $[\xi_i,\xi_k]=0$. Moreover $\xi_i1=0$ and $B=\C[F]\cdot1$. (ii) We have $[\sigma f,h]_\sigma=\sigma\br{\sigma f_1,\ldots,\sigma f_{n-1},h}=\sigma(X_{\sigma f}h)$, and $X_{\sigma f}=F^!X_f$ by Proposition~\ref{prop:lift}.
\end{proof}

For the shear of Example~\ref{ex:shear} we have $(dF)^{-1}=\begin{pmatrix}1&-2y\\0&1\end{pmatrix}$, so $\xi_1=\partial_x$ and $\xi_2=\partial_y-2y\partial_x$. Indeed $\xi_2(x+y^2)=0$ and $\xi_2y=1$.

\begin{proposition}\label{prop:Q}
Let $F\colon X\to Y$ be a dominant local isomorphism of smooth affine varieties and $Q=A/B$. For $y\in Y$,
\[ Q\otimes_Bk(y)\cong\O\big(F^{-1}(y)\big)/\C\cdot1, \]
which has dimension $\#F^{-1}(y)-1$ if $F^{-1}(y)\ne\emptyset$ and $0$ otherwise. Moreover $Q=0$ if and only if $F$ is an isomorphism.
\end{proposition}
\begin{proof}
Tensoring $B\to A\to Q\to0$ with $k(y)$ over $B$ is right exact and gives $Q\otimes k(y)=\operatorname{coker}(k(y)\to A\otimes_Bk(y))$. Here $A\otimes_Bk(y)=\O(F^{-1}(y))$ is the algebra of functions on a finite reduced set by (G1), and the map is the inclusion of constants. If $Q=0$ then $F^*$ is surjective, so $F$ is a closed immersion. It is also open and dominant, and $Y$ is irreducible, so $F$ is an isomorphism. The converse is clear.
\end{proof}

\begin{example}\label{ex:toy}
Let $X=\Spec\C[t,t^{-1},(1+t)^{-1}]$, $Y=\mathbb G_m=\Spec\C[s,s^{-1}]$ and $F(t)=t^2$. This is a local isomorphism of degree two. The fibre over $s\ne1$ is $\{\pm\sqrt s\}$, but the fibre over $s=1$ is $\{1\}$ because the point $t=-1$ has been removed. So $F$ is not proper over $s=1$. By Proposition~\ref{prop:Q} the fibre of $Q$ has dimension one for $s\ne1$ and zero at $s=1$. The connection is $\nabla_{\partial_s}=\frac1{2t}\partial_t$, and going once around $s=0$ interchanges the two sheets, so $Q$ over $Y\setminus\{1\}$ is the rank one local system with monodromy $-1$.
\end{example}

\subsection{Non-properness and monodromy}

\begin{proposition}\label{prop:nonproper}
A Keller map $F\colon\A^n\to\A^n$ that is not an automorphism is not proper.
\end{proposition}
\begin{proof}
Suppose $F$ is proper. It is quasi-finite by (G1), hence finite by (G2), hence proper in the classical topology. As a local homeomorphism it is then a covering of $\C^n$ with connected source by (G3), hence a homeomorphism because $\C^n$ is simply connected. Thus $F$ is injective, hence an automorphism by (G4).
\end{proof}

Let $S_F\subset Y$ be the set of points $y$ that have no neighbourhood $V$ with $F^{-1}(V)\to V$ proper; we call it the non-properness set. It is closed, and $S_F=\emptyset$ if and only if $F$ is proper. By Proposition~\ref{prop:nonproper}, $S_F\ne\emptyset$ for a Keller map that is not an automorphism. Let $d=\deg F$ and let $X\subset X^\nu\to Y$ be the factorisation of (G5). By (G5), $\#F^{-1}(y)\le d$ for all $y$. If $\#F^{-1}(y)=d$, then $F^{-1}(y)$ is the whole fibre of $X^\nu$ over $y$, so the closed set $X^\nu\setminus X$ has image missing $y$, and $F$ is finite near $y$. Conversely, let $V$ be a Zariski open neighbourhood of $y$ over which $F$ is finite. It is connected because $Y$ is irreducible, and $F$ is a finite \'etale covering over it, so all its fibres have the same cardinality, which is $d$. Therefore
\begin{equation}\label{eq:SF}
S_F=\{y\in Y:\ \#F^{-1}(y)<d\}.
\end{equation}
The same argument shows $S_F=\pi(X^\nu\setminus X)$, where $\pi\colon X^\nu\to Y$ is the finite map of (G5).

\begin{proposition}\label{prop:mono}
Let $U=Y\setminus S_F$. Then $F_*\O_X|_U$ is a flat vector bundle of rank $d$ whose horizontal sections form the permutation local system of the $d$-sheeted covering $F^{-1}(U)\to U$. The bundle $Q|_U$ is the flat bundle of rank $d-1$ associated to this permutation representation of $\pi_1(U)$ modulo the trivial subrepresentation.
\end{proposition}
\begin{proof}
Over $U$ the map $F$ is a finite local isomorphism, hence a $d$-sheeted covering map. A section $h$ over a small disc $V\subset U$ is a function on $F^{-1}(V)=\bigsqcup_{i=1}^dV_i$, and $\nabla h=0$ if and only if $h$ is locally constant, that is, a $d$-tuple of constants. The subbundle $B|_U=\O_U\cdot1$ is the trivial subrepresentation, and the quotient is $Q|_U$.
\end{proof}

\section{Ramification at infinity}\label{sec:infinity}

Let $F$ be a Keller map. A projective extension is a pair $(\bar X,\bar F)$ with $\bar X$ a smooth projective variety containing $\A^n$ as a dense open subset and $\bar F\colon\bar X\to\P^n$ a morphism with $\bar F|_{\A^n}=F$. Write $D_\infty=\bar X\setminus\A^n=\bigcup_iD_i$ for the boundary. Projective extensions exist: take a smooth compactification and resolve the indeterminacy of $F$, which lies in the boundary, by blow-ups with centres in the boundary. Let $H_\infty=\P^n\setminus\A^n$. The form $\omega=dx_1\wedge\cdots\wedge dx_n$ is a rational section of $\Omega^n_{\bar X}$ without zeros or poles on $\A^n$, so $K_{\bar X}=\divv\omega=-\sum_ia_iD_i$ with $a_i=-\ord_{D_i}\omega$. Likewise $K_{\P^n}=\divv(du_1\wedge\cdots\wedge du_n)=-(n+1)H_\infty$.

\begin{theorem}\label{thm:C}
Let $F$ be a Keller map and $(\bar X,\bar F)$ a projective extension. The ramification divisor of $\bar F$ is
\[ R^\infty_F:=K_{\bar X}-\bar F^*K_{\P^n}=(n+1)\bar F^*H_\infty-\sum_ia_iD_i\ \ge0, \]
and it is supported on $D_\infty$.
\end{theorem}
\begin{proof}
The map $\bar F$ is generically finite, so (G6) gives $R_{\bar F}=K_{\bar X}-\bar F^*K_{\P^n}\ge0$. As rational forms $\bar F^*(du_1\wedge\cdots\wedge du_n)=c\,\omega$. So on $\A^n$ the Jacobian section is the non-zero constant $c$ in the frames $du_1\wedge\cdots\wedge du_n$ and $\omega$, and $R_{\bar F}$ misses $\A^n$. The second expression follows by inserting the canonical divisors.
\end{proof}

\begin{theorem}\label{thm:nonvanish}
If $R^\infty_F=0$ then $\bar F$ is an isomorphism and $F$ is an automorphism. Hence for a Keller map $F$ that is not an automorphism, $R^\infty_F\ne0$ for every projective extension.
\end{theorem}
\begin{proof}
The condition $R^\infty_F=0$ means that $\bar F^*\Omega^n_{\P^n}\to\Omega^n_{\bar X}$ is an isomorphism, so $\bar F$ is a local isomorphism by (G1). It is proper, hence finite by (G2), and a covering of $\P^n$ by (G3). Since $\P^n$ is simply connected and $\bar X$ is connected, $\bar F$ is a homeomorphism, and being a finite morphism that is bijective and \'etale it is an isomorphism. Then $F=\bar F|_{\A^n}$ is injective, hence an automorphism by (G4).
\end{proof}

\begin{corollary}\label{cor:Jel}
Let $F$ be a Keller map that is not an automorphism, $(\bar X,\bar F)$ a projective extension, and $X'=\bar F^{-1}(\A^n)$. Then $S_F=\bar F(D_\infty)\cap\A^n$ is non-empty and closed of pure codimension one~\cite{Jelonek}. Let $\Sigma$ be the set of boundary components $D_i$ with $\bar F(D_i)\not\subset H_\infty$. Then $\Sigma\ne\emptyset$ and $S_F=\bigcup_{D_i\in\Sigma}\bar F(D_i)\cap\A^n$, and the coefficient of $D_i\in\Sigma$ in $R^\infty_F$ is $\ord_{D_i}\omega=-a_i\ge0$.
\end{corollary}
\begin{proof}
The map $F'=\bar F|_{X'}\colon X'\to\A^n$ is proper, being a base change of a proper morphism, and $F=F'|_{\A^n}$. The set $F'(X'\setminus\A^n)=\bar F(D_\infty)\cap\A^n$ is closed. If $y$ lies outside it, then $F^{-1}(V)=F'^{-1}(V)$ for a neighbourhood $V$ of $y$, so $F$ is proper over $V$. If $y=F'(p)$ with $p\in X'\setminus\A^n$ and $V\ni y$ is open, then $F^{-1}(V)$ is a dense proper open subset of the irreducible variety $F'^{-1}(V)$, which contains $p$. If $F^{-1}(V)\to V$ were proper, then $F^{-1}(V)\to F'^{-1}(V)$ would be proper because $F'^{-1}(V)\to V$ is separated. A proper open immersion is also closed, hence an equality on the irreducible target, which is false. So $S_F=\bar F(D_\infty)\cap\A^n$. By Proposition~\ref{prop:nonproper} this set is non-empty, hence $\Sigma\ne\emptyset$. For $D_i\in\Sigma$ we have $\ord_{D_i}\bar F^*H_\infty=0$, so Theorem~\ref{thm:C} gives the coefficient $-a_i$. Purity of codimension one is Jelonek's theorem.
\end{proof}

\begin{remark}\label{rem:BGV}
Let $X_F$ be the normalization of the target in the function field of the source, with its open embedding $\A^n\hookrightarrow X_F$ and finite map $\pi\colon X_F\to\A^n$~\cite{BGV}. By~\eqref{eq:SF} and (G5), $S_F=\pi(X_F\setminus\A^n)$. A projective extension in which the boundary divisors lying over $X_F\setminus\A^n$ are realised dominates the pair $(X_F,\pi)$, and Corollary~\ref{cor:Jel} is the statement for such a model. Components of $X_F\setminus\A^n$ in the sense of~\cite{BGV} therefore correspond to the elements of $\Sigma$ that are not contracted by the passage from $\bar X$ to $X_F$.
\end{remark}

\begin{remark}\label{rem:Pn}
If $\bar X=\P^n$ and $\bar F\colon\P^n\to\P^n$ extends a polynomial map $F$, then $\bar F^{-1}(\A^n)=\A^n$ because $\bar F$ is finite and $\bar F^{-1}(H_\infty)$ is a hypersurface containing $H_\infty$, hence equal to it. So $F$ is proper, and a Keller map $F$ is then an automorphism by Proposition~\ref{prop:nonproper}. Thus a Keller map that is not an automorphism never extends to a morphism $\P^n\to\P^n$. Every projective extension needs blow-ups, and on them $R^\infty_F\ne0$.
\end{remark}

\begin{example}\label{ex:toy2}
In Example~\ref{ex:toy} take $\bar X=\P^1\supset X$, whose boundary is $\{0,\infty,-1\}$, and $\bar F(t)=t^2$. The ramification divisor is $[0]+[\infty]$, as Riemann--Hurwitz gives $-2=2\cdot(-2)+2$, and it is supported on the boundary. The non-properness set is $\bar F(\{0,\infty,-1\})\cap Y=\{1\}$, as computed directly.
\end{example}

\section{Counterexamples in dimension at least three}\label{sec:alpoge}

\subsection{What holds for every counterexample}

A counterexample to the Jacobian conjecture is a non-injective Keller map $F\colon\A^n\to\A^n$; by (G4) this is the same as a Keller map that is not an automorphism. Counterexamples exist for every $n\ge3$~\cite{Gao}. None exists for $n=1$, and the case $n=2$ is open. Everything proved so far applies to each of them.

\begin{proposition}\label{prop:general}
Let $F\colon\A^n\to\A^n$ be a counterexample, $Y=\A^n$ the target, $d=\deg F$ and $U=Y\setminus S_F$.
\begin{enumerate}
\item[(i)] $2\le d<\infty$, and the twist of $F$ is essential.
\item[(ii)] The set $S_F=\{y:\#F^{-1}(y)<d\}$ is non-empty and closed of pure codimension one. Every projective extension has $R^\infty_F\ne0$ and $S_F=\bar F(D_\infty)\cap\A^n$.
\item[(iii)] The map $F^{-1}(U)\to U$ is a connected $d$-sheeted covering, its monodromy group $G\subset S_d$ is transitive, and $Q|_U$ is the flat bundle of rank $d-1$ attached to the permutation representation of $G$ modulo the trivial one. In particular $Q\otimes k(y)\ne0$ for every $y$ in the image of $F$ with $\#F^{-1}(y)\ge2$.
\item[(iv)] If $\alpha,\beta$ are automorphisms of $\A^n$ with constant Jacobian, then $\alpha\circ F\circ\beta$ is again a counterexample with $S_{\alpha F\beta}=\alpha(S_F)$ and the same $d$ and $G$. The map $F\times\mathrm{id}_{\A^k}\colon\A^{n+k}\to\A^{n+k}$ is a counterexample with $S=S_F\times\A^k$ and the same $d$ and $G$.
\end{enumerate}
\end{proposition}
\begin{proof}
(i) The degree is finite because $F$ is quasi-finite and dominant. If $d=1$ then $F$ would be a quasi-finite birational morphism onto a normal variety, hence an open immersion by (G5), hence injective, contrary to assumption. The twist is essential by Proposition~\ref{prop:untwist}. (ii) This is Proposition~\ref{prop:nonproper}, formula~\eqref{eq:SF} and Corollary~\ref{cor:Jel}. (iii) The covering is connected because $F^{-1}(U)$ is a non-empty open subset of the irreducible variety $\A^n$, so $G$ is transitive. The rest is Propositions~\ref{prop:mono} and~\ref{prop:Q}. (iv) Automorphisms preserve properness and the cardinality of fibres, and they multiply the Jacobian by a constant. For the product, $(F\times\mathrm{id})^{-1}(V\times W)=F^{-1}(V)\times W\to V\times W$ is proper if and only if $F^{-1}(V)\to V$ is, by base change and restriction to a point of $W$. So $S_{F\times\mathrm{id}}=S_F\times\A^k$, while fibres and degree are unchanged.
\end{proof}

\subsection{Cubic forms with a marked root}\label{sec:cubic}

From now on $F$ is the map~\eqref{eq:Alpoge}. That it is a Keller map with $\det JF=-2$, and that it identifies the three points $p_1=(0,0,-\tfrac14)$, $p_2=(1,-\tfrac32,\tfrac{13}2)$, $p_3=(-1,\tfrac32,\tfrac{13}2)$ with $q=(-\tfrac14,0,0)$, is established in~\cite{Alpoge,RHQ} and not repeated. We write $A=\C[x,y,z]$, $B=\C[F_1,F_2,F_3]\subset A$, $K=\C(x,y,z)\supset L=\C(F_1,F_2,F_3)$, and $u=(u_1,u_2,u_3)$ for coordinates on the target $Y=\A^3$.

For $u\in Y$ let
\[ p_u(\mu,\eta)=u_3\mu^3-2\mu^2\eta+u_2\mu\eta^2-2u_1\eta^3 \]
be a binary cubic form. It is never zero because of the coefficient $-2$. Let
\[ \CC=\big\{(u,\rho)\in Y\times\P^1:\ \rho\text{ is a simple root of }p_u\big\}. \]
In the notation of~\cite{BGVai}, $(p,q,r)=(u_1,u_2,u_3)$, the form $\Phi_{p,q,r}(S,T)$ equals $-p_u(T,S)$, the root $[S:T]=[x:1+xy]$ corresponds to $\rho=(1+xy:x)$, the incidence variety $X=V(\Phi)$ is the closure of $\CC$ obtained by allowing multiple roots, and $\Delta$ below is their $\Delta$. We write $\pi\colon X\to Y$ for the finite flat cubic cover.

\begin{theorem}[{\cite[\S4]{Alpoge}, \cite[Prop.\,2.1, Prop.\,2.2, Prop.\,3.1]{BGVai}}]\label{thm:cubic}
The map $\Psi\colon\A^3\to\CC$, $\Psi(x,y,z)=\big(F(x,y,z),(1+xy:x)\big)$, is an isomorphism of varieties, and $F=\mathrm{pr}_Y\circ\Psi$. The projection $\pi\colon X\to Y$ is finite flat of degree three with $X$ smooth and affine, $\CC\subset X$ is the open locus of simple roots, and $X$ is the normalization $X_F$ of $Y$ in $K$ in the sense of~\cite{BGV}. The boundary $D=X\setminus\CC$ is the locus of multiple roots, it is isomorphic to $\A^2$, and it is the unique irreducible ramification divisor of $\pi$.
\end{theorem}

We use the explicit form of $\Psi$ on the chart $\{x\ne0\}$ in what follows. Put $w=1+xy$ and, on $\{x\ne0\}$, $m=w/x=y+1/x$, and $P(\mu)=F_3\mu^3-2\mu^2+F_2\mu-2F_1$.

\begin{lemma}\label{lem:chart}
On $\{x\neq0\}$ we have
\begin{equation}\label{eq:cubic}
F_2=y+3F_1/m,\qquad P(m)=0,\qquad x=\frac{2}{P'(m)}.
\end{equation}
Let $T=\A^3$ with coordinates $(\beta,\gamma,\mu)$, $\theta=3\gamma\mu^2-4\mu+\beta$, and $\psi=(F_2,F_3,m)$. Then $\psi$ is an isomorphism of $X_0=\{x\ne0\}$ onto $T\setminus\{\theta=0\}$ with $\theta\circ\psi=2/x$, and its inverse is
\[ \Theta(\beta,\gamma,\mu)=\Big(\frac2\theta,\ \mu-\frac\theta2,\ \frac{5\theta^2}4-\frac{3\mu\theta}2-\frac{\gamma\theta^3}8\Big). \]
Under $\Psi$ the set $X_0$ corresponds to the part of $\CC$ with $\rho\neq\infty$, and $T=\{(\beta,\gamma,\mu)\}$ is the chart of $X$ in which $D=\{\theta=0\}$.
\end{lemma}
\begin{proof}
Since $4+3xy=1+3w$ and $2-3xy=5-3w$, we have $F_1=w[w^2z+y^2(1+3w)]$, $F_2=y+3x[w^2z+y^2(1+3w)]$ and $F_3=x(5-3w)-x^3z$. The first two give $F_2=y+3F_1/m$. Solving the third for $z$ and inserting in $F_1$, with $w^3/x^3=m^3$ and $(w-1)^2(1+3w)+w^2(5-3w)=w+1$, gives $F_1+m^3F_3=w(w+1)/x^2=2m^2-my$, which with $my=mF_2-3F_1$ is $P(m)=0$. From $1/x=m-y=(m^2-F_2m+3F_1)/m$ and $3F_1=\tfrac32(F_3m^3-2m^2+F_2m)$ one gets $m^2-F_2m+3F_1=\tfrac m2P'(m)$, which is the third relation. Hence $\theta\circ\psi=P'(m)=2/x$. The formulas for $\Theta$ then read $x=2/\theta$, $y=m-1/x$ and $z=(x(5-3w)-F_3)/x^3$ with $w=mx$, so $\Theta\circ\psi=\mathrm{id}$; conversely $m\circ\Theta=\mu$, $F_3\circ\Theta=\gamma$ by the definition of $z$, and eliminating $f_1=F_1\circ\Theta$ from $2f_1=\gamma\mu^3-2\mu^2+f_2\mu$ and $\mu f_2=\mu^2-\theta\mu/2+3f_1$ gives $f_2=4\mu+\theta-3\gamma\mu^2=\beta$. These are identities of regular functions on an irreducible variety, so division by $\mu$ is harmless. A multiple root $\mu$ of $p_u(\mu,1)$ is a point with $P(\mu)=P'(\mu)=0$, that is with $\theta=0$, and a multiple root is never at $\infty$ because $\partial_\eta p_u(1,0)=-2$.
\end{proof}

\begin{remark}\label{rem:keller-explained}
In the chart $\CC_0\cong T\setminus\{\theta=0\}$ the map $F$ is $G(\beta,\gamma,\mu)=\big(\tfrac12(\gamma\mu^3-2\mu^2+\beta\mu),\beta,\gamma\big)$ with $\det dG=\partial_\mu u_1=\theta/2$ and $\psi^*\theta=2/x$. So $J_F=\tfrac12(\theta\circ\psi)J_\psi=\tfrac1xJ_\psi$, and the constancy $J_F=-2$ says exactly that $J_\psi=-2x$. The Keller property is thus the compatibility of the volume form on the space of cubic forms with the volume form on $\CC\cong\A^3$, namely $\omega=-\tfrac12F^*(du_1\wedge du_2\wedge du_3)$. The plane $\{x=0\}$ is the locus of roots at $\infty$.
\end{remark}

\subsection{Degree, fibres, non-properness and monodromy}

Let $\Delta(u)=u_2^2-u_2^3u_3-16u_1-27u_1^2u_3^2+18u_1u_2u_3$, so that $\Disc(p_u)=4\Delta(u)$, let $\Sigma=\{\Delta=0\}$, and let $\Gamma\subset\Sigma$ be the curve of triple-root forms, $\Gamma=\{12u_1=u_2^2,\ 3u_2u_3=4\}\cong\mathbb G_m$.

\begin{theorem}\label{thm:fibres}
\begin{enumerate}
\item[(i)]{\rm \cite[Thm.\,4.5]{BGVai}} The fibre $F^{-1}(u)$ is in bijection with the set of simple roots of $p_u$ in $\P^1$. Hence $\#F^{-1}(u)=3,1,0$ according as $u\notin\Sigma$, $u\in\Sigma\setminus\Gamma$, $u\in\Gamma$. In particular $F(\A^3)=\A^3\setminus\Gamma$, and the non-properness set is the irreducible discriminant surface $S_F=\Sigma$.
\item[(ii)] $K=L(m)$ with $m=y+1/x$, the minimal polynomial of $m$ over $L$ is $P(\mu)=F_3\mu^3-2\mu^2+F_2\mu-2F_1$, and $\deg F=[K:L]=3$.
\item[(iii)] The Galois closure of $K/L$ has group $S_3$~\cite[Prop.\,8.1]{BGVai}, and $\Aut(K/L)=1$, so $F$ has no non-trivial deck transformation although it has degree three. The monodromy representation $\pi_1(Y\setminus\Sigma)\to S_3$ is surjective. The bundle $Q|_{Y\setminus\Sigma}$ is the flat bundle of the standard two-dimensional irreducible representation of $S_3$, and $Q\otimes_Bk(u)$ has dimension $2$ for $u\notin\Sigma$ and dimension $0$ for $u\in\Sigma$.
\item[(iv)] At $q=(-\tfrac14,0,0)$ we have $\Delta(q)=4$ and $p_q=\eta(-2\mu^2+\tfrac12\eta^2)$, whose simple roots are $\infty$ and $\mu=\pm\tfrac12$. They correspond to $p_1$, where $x=0$, and to $p_2,p_3$, where $m=\mp\tfrac12$ and $x=-1/(2m)=\pm1$. Moreover $Q\otimes_Bk(q)\cong\C^3/\C(1,1,1)$, with basis the classes of $x$ and $z$.
\end{enumerate}
\end{theorem}
\begin{proof}
(i) The stratification of the fibres is~\cite[Thm.\,4.5]{BGVai}; it follows from Theorem~\ref{thm:cubic} because a binary cubic form has three simple roots if $\Disc\ne0$, a double and a simple root if $\Disc=0$ and the form is not a cube, and a triple root exactly on $\Gamma$, where $p_u=u_3(\mu-a\eta)^3$ forces $3u_3a=2$, $u_2=2a$, $u_1=a^2/3$. By~\eqref{eq:SF} with $d=3$, $S_F=\{\#F^{-1}<3\}=\Sigma$. Irreducibility of $\Sigma$ is~\cite[Prop.\,4.2]{BGVai}: $\Sigma=\pi(D)$ with $D\cong\A^2$.

(ii) By~\eqref{eq:cubic}, $m$ is a root of the stated cubic over $L$. By Lemma~\ref{lem:chart}, $x,y,z\in L(m)$, being the coordinates of $\Theta(F_2,F_3,m)$. So $K=L(m)$ and $[K:L]\le3$. By (i) there are $u$ with $\#F^{-1}(u)=3$, so $\deg F\ge3$ by (G5). Hence $\deg F=3$ and the cubic is irreducible.

(iii) The coefficients of $m^3-\tfrac{2}{F_3}m^2+\tfrac{F_2}{F_3}m-\tfrac{2F_1}{F_3}$ are $e_1=2/F_3$, $e_2=F_2/F_3$ and $e_3=2F_1/F_3$, and $L=\C(e_1,e_2,e_3)$ since $F_3=2/e_1$, $F_2=2e_2/e_1$ and $F_1=e_3/e_1$. The $e_i$ are algebraically independent because $F_1,F_2,F_3$ are, so $m$ is a root of the generic cubic, whose Galois group is $S_3$~\cite[Ch.\,VI, \S2]{Lang}. The monodromy group of the connected covering $F^{-1}(Y\setminus\Sigma)\to Y\setminus\Sigma$ is this Galois group~\cite[\S I]{Harris}. For a cubic extension with Galois group $S_3$ we have $\Aut(K/L)=N_{S_3}(S_2)/S_2=1$. The statement on $Q$ over $Y\setminus\Sigma$ is Proposition~\ref{prop:mono}, because the permutation representation of $S_3$ on three letters is the sum of the trivial and the standard representation. Over $\Sigma$ the fibre has at most one point by (i), so $Q\otimes k(u)=0$ by Proposition~\ref{prop:Q}.

(iv) Evaluate $\Delta$ and $p_q$ directly. For the point $p_2$ we have $m=-\tfrac32+1=-\tfrac12$, and for $p_3$ we have $m=\tfrac32-1=\tfrac12$. Here $P=-2\mu^2+\tfrac12$, so $x=2/P'(m)=-1/(2m)$. By Proposition~\ref{prop:Q} the fibre of $Q$ at $q$ is the quotient of the functions on the three points by constants. The values $(x(p_i))=(0,1,-1)$ and $(z(p_i))=(-\tfrac14,\tfrac{13}2,\tfrac{13}2)\equiv(-\tfrac{27}4,0,0)$ modulo $(1,1,1)$ are linearly independent in $\C^3/\C(1,1,1)$.
\end{proof}

The dimension $\dim Q\otimes k(q)=2$ is thus the generic rank of $Q$ and nothing special at $q$. What is special is the global structure. Over the complement of the discriminant surface $\Sigma$, which has a cuspidal edge along $\Gamma$~\cite[Thm.\,4.4]{BGVai}, we have a flat $S_3$-bundle. Over $\Sigma\setminus\Gamma$ the fibre of $F$ is a single point, because two of the three marked roots merge into a double root and leave $\A^3$, and $Q$ is zero there; over $\Gamma$ the fibre is empty. So $Q$ has rank $2$ on the open set $\A^3\setminus\Sigma$ and rank $0$ on its complement, and its support is the closure of the locus where $F$ is a three-sheeted covering.

\subsection{The twisted ternary algebra}

Take $\kappa=i/\sqrt2$, so that $\kappa^2\cdot(-2)=1$, and put $\sigma=\kappa F^*$. Write $u_i$ for the coordinate functions of the target, so that $\sigma(u_i)=\kappa F_i$. By Proposition~\ref{prop:hom}, $\mathcal A=(A,[\cdot,\cdot,\cdot]_\sigma,\sigma)$ is a multiplicative Hom-Nambu algebra with $[f,g,h]_\sigma=\kappa\,\br{f,g,h}\circ F$. For instance $[x,y,z]_\sigma=\kappa$ and $[\sigma u_1,\sigma u_2,\sigma u_3]_\sigma=-2\kappa^4=-\tfrac12=\sigma[x,y,z]_\sigma$, as multiplicativity requires. With $\theta=3\gamma\mu^2-4\mu+\beta$ as before, $A[1/x]=\C[\beta,\gamma,\mu][1/\theta]$ with $\beta=F_2$, $\gamma=F_3$, $\mu=m$, and
\[ x=\frac2\theta,\qquad y=\mu-\frac\theta2,\qquad z=\tfrac54\theta^2-\tfrac32\mu\theta-\tfrac18\gamma\theta^3. \]

\begin{theorem}\label{thm:twisted}
\begin{enumerate}
\item[(i)] The lifts $\xi_i=F^!\partial_{u_i}$ are the columns of $-\tfrac12\adj(dF)$, hence polynomial vector fields on $\A^3$. In the coordinates $(\beta,\gamma,\mu)$ on $A[1/x]$,
\[ \xi_1=\frac2\theta\,\partial_\mu,\qquad \xi_2=\partial_\beta-\frac\mu\theta\,\partial_\mu,\qquad \xi_3=\partial_\gamma-\frac{\mu^3}\theta\,\partial_\mu, \]
that is, $\xi_i\mu=-\partial_{u_i}P/P'(\mu)$, so the $\xi_i$ move the marked root along the family of cubic forms. Moreover
\[ [\sigma u_1,\sigma u_2,h]_\sigma=\sigma(\xi_3h),\quad[\sigma u_2,\sigma u_3,h]_\sigma=\sigma(\xi_1h),\quad[\sigma u_3,\sigma u_1,h]_\sigma=\sigma(\xi_2h). \]
\item[(ii)] On $A[1/x]$ the Jacobian brackets in the coordinates $(x,y,z)$ and $(\beta,\gamma,\mu)$ are related by $\br{\cdot}_{xyz}=-2x\,\br{\cdot}_{\beta\gamma\mu}=-\tfrac4\theta\br{\cdot}_{\beta\gamma\mu}$, and $G^*$ is conformal with factor $\theta/2$.
\item[(iii)] The map $\sigma\colon\mathcal A\to\mathcal A$ is an injective morphism of Hom-Nambu algebras with image $(B,[\cdot]_\sigma|_B,\sigma|_B)\cong\mathcal A$. Hence $\mathcal A\supset\sigma\mathcal A\supset\sigma^2\mathcal A\supset\cdots$, each term a copy of $\mathcal A$. On fields of fractions $K\supset L\supset\sigma(L)\supset\cdots$ with $[K:\sigma^k(K)]=3^k$. The extension $B\subset A$ is the cubic extension $L(m)/L$, localised and glued as in Theorem~\ref{thm:cubic}: $A[1/x]=(B[\mu]/P(\mu))[1/P'(\mu)]$, and $A=A[1/x]\cap A[1/w]$, where $A[1/w]=(B[\eta]/p(1,\eta))[1/\partial_\eta p(1,\eta)]$ is the ring of the second chart.
\end{enumerate}
\end{theorem}
\begin{proof}
(i) The fields $\xi_i$ are the columns of $(dF)^{-1}=-\tfrac12\adj(dF)$, a polynomial matrix. On $A[1/x]$ we have $F=(g,\beta,\gamma)$ with $g=\tfrac12(\gamma\mu^3-2\mu^2+\beta\mu)$, and $\xi_i$ is characterised by $\xi_i\beta=\delta_{i2}$, $\xi_i\gamma=\delta_{i3}$ and $\xi_ig=\delta_{i1}$. The first two conditions give $\xi_i=a_i\partial_\mu+\delta_{i2}\partial_\beta+\delta_{i3}\partial_\gamma$. The third condition reads $a_ig_\mu+\delta_{i2}g_\beta+\delta_{i3}g_\gamma=\delta_{i1}$ with $g_\mu=\theta/2$, $g_\beta=\mu/2$ and $g_\gamma=\mu^3/2$, which gives $a_1=2/\theta$, $a_2=-\mu/\theta$ and $a_3=-\mu^3/\theta$. Writing $P(\mu)=u_3\mu^3-2\mu^2+u_2\mu-2u_1$ we have $\partial_{u_1}P=-2$, $\partial_{u_2}P=\mu$, $\partial_{u_3}P=\mu^3$ and $P'=\theta$, so $a_i=-\partial_{u_i}P/P'$. The three bracket identities are Proposition~\ref{prop:conn}(ii), because $X_{u_1,u_2}=\partial_{u_3}$, $X_{u_2,u_3}=\partial_{u_1}$ and $X_{u_3,u_1}=\partial_{u_2}$.

(ii) By Remark~\ref{rem:keller-explained}, $\psi^*(d\beta\wedge d\gamma\wedge d\mu)=J_\psi\,dx\wedge dy\wedge dz=-2x\,dx\wedge dy\wedge dz$, which gives the first factor. Since $\psi^*\theta=2/x$ it equals $-4/\theta$. The factor $\theta/2=\partial_\mu u_1$ for $G^*$ is Lemma~\ref{lem:conf}.

(iii) We have $\sigma[f,g,h]_\sigma=\sigma^2\br{f,g,h}=\sigma\br{\sigma f,\sigma g,\sigma h}=[\sigma f,\sigma g,\sigma h]_\sigma$, and $\sigma\sigma=\sigma\sigma$, so $\sigma$ is a morphism. It is injective with image $B$ by Proposition~\ref{prop:untwist}, so it is an isomorphism onto $(B,[\cdot]_\sigma|_B,\sigma|_B)$. The degrees are those of Theorem~\ref{thm:fibres}(ii), multiplied along the tower. The description of $A[1/x]$ is Lemma~\ref{lem:chart}. Since $w=1+xy$ and $x$ generate the unit ideal of $A$, the open sets $\{x\ne0\}$ and $\{w\ne0\}$ cover $\A^3$ and $A=A[1/x]\cap A[1/w]$. The chart $\{w\ne0\}$ corresponds under $\Psi$ to the second chart $\{\rho\ne0\}$ of $\CC$, which gives the second ring.
\end{proof}

The coefficients $a_i$ in (i) have the common denominator $\theta$, and the numerators do not vanish along $\{\theta=0\}$ (for $\xi_1$ the numerator is the constant $2$). This gives the following description of the connection of Proposition~\ref{prop:conn} on the finite completion.

\begin{proposition}\label{prop:polar}
Regard the fields $\xi_1,\xi_2,\xi_3$ as meromorphic vector fields on the cubic completion $X$ of Theorem~\ref{thm:cubic}. They are regular on $\CC=\A^3$ and have a simple pole along the ramification divisor $D=\{\theta=0\}$ of $\pi\colon X\to Y$, and no other pole. Equivalently, $\theta\,\xi_i$ is regular on $X$ and does not vanish identically on $D$ for each $i$, so the flat connection $\nabla$ on $\pi_*\O_X|_{Y\setminus\Sigma}$ has a pole along $\Sigma=\pi(D)$ with the local type of the ramification of $\pi$ at $D$.
\end{proposition}
\begin{proof}
The chart $T$ of Lemma~\ref{lem:chart} contains $D$, since a multiple root is never at $\infty$. By Theorem~\ref{thm:twisted}(i) the coefficients of the $\xi_i$ in the frame $\partial_\beta,\partial_\gamma,\partial_\mu$ of $T$ are polynomials in $\beta,\gamma,\mu$ divided by $\theta$, and $\theta\,\xi_1=2\partial_\mu$, $\theta\,\xi_2=\theta\partial_\beta-\mu\partial_\mu$, $\theta\,\xi_3=\theta\partial_\gamma-\mu^3\partial_\mu$. On $D=\{\theta=0\}$ these restrict to $2\partial_\mu$, $-\mu\partial_\mu$ and $-\mu^3\partial_\mu$. The first is nowhere zero. The function $\mu$ is not identically zero on the irreducible surface $D$, because $\mu=0$ and $\theta=0$ give $\beta=0$, which is a proper closed subset of $D$ ($D$ is a graph over the $(\gamma,\mu)$-plane). So each $\theta\xi_i$ is regular and non-zero along $D$. Off $D$ the fields are regular because $\theta$ is a unit there.
\end{proof}

\begin{remark}\label{rem:K_X}
Under the identification of Theorem~\ref{thm:cubic} one finds from Theorem~\ref{thm:twisted}(ii) that $\omega=-\tfrac\theta4\,d\beta\wedge d\gamma\wedge d\mu$ on $T\setminus\{\theta=0\}$. This is the expression of the Riemann--Hurwitz relation $K_X=\pi^*K_Y+D$ of~\cite[\S6]{BGVai} in the frame $\omega$ of $\A^3$: the form $\omega$ vanishes to order one along $D$ relative to a generator of $\Omega^3_X$, and $\pi^*(du_1\wedge du_2\wedge du_3)=-2\omega$ has the corresponding order one pole.
\end{remark}

\begin{remark}\label{rem:higher}
For $n>3$ the map $F_n=F\times\mathrm{id}_{\A^{n-3}}$ is a Keller map with $J=-2$. Theorems~\ref{thm:cubic}--\ref{thm:twisted} hold for it with $\CC\times\A^{n-3}$ in place of $\CC$, with $S_{F_n}=\Sigma\times\A^{n-3}$, and the twist is normalised by $\kappa^{n-1}\cdot(-2)=1$.
\end{remark}

\subsection{The boundary}

By Corollary~\ref{cor:Jel} the non-properness set $S_F=\Sigma$ is the image of boundary components of any projective extension. The structure above identifies the relevant component. Let $\nu$ be the divisorial valuation of $K$ given by the order of vanishing along the surface $\{\theta=0\}\subset T$, where $T\supset T\setminus\{\theta=0\}\cong X_0$ via $\psi$. Then $\nu(x)=\nu(2/\theta)=-1<0$, so the centre of $\nu$ on every model containing $\A^3$ as an open subset lies in the boundary. Its centre on the finite completion $X$ is the generic point of $D$.

\begin{proposition}\label{prop:nu}
Let $(\bar X,\bar F)$ be a projective extension in which $\nu$ is realised by a boundary component $D_\nu$. Such extensions exist: blow up the centre of $\nu$ repeatedly, resolve, and keep $D_\nu$. Then $\bar F(D_\nu)\cap\A^3=S_F=\Sigma$, and the coefficient of $D_\nu$ in $R^\infty_F$ is exactly $1$.
\end{proposition}
\begin{proof}
The centre of $\nu$ on $Y$ is the closure of $G(\{\theta=0\})=\Sigma$ by Theorem~\ref{thm:fibres}(i). Hence $D_\nu\in\Sigma$ and $\bar F(D_\nu)\cap\A^3=S_F$. On $T\setminus\{\theta=0\}$ we have $\omega=dx\wedge dy\wedge dz=-\tfrac1{2x}\,d\beta\wedge d\gamma\wedge d\mu=-\tfrac\theta4\,d\beta\wedge d\gamma\wedge d\mu$ by Theorem~\ref{thm:twisted}(ii). The form $d\beta\wedge d\gamma\wedge d\mu$ generates $\Omega^3_T$ at the generic point of $\{\theta=0\}$, where $\theta$ vanishes to first order because it is linear in $\beta$. The order of a rational top form along a divisorial valuation depends only on the valuation ring, so $\ord_{D_\nu}\omega=\nu(\theta)=1$. Corollary~\ref{cor:Jel} gives the coefficient.
\end{proof}

Thus the ramification divisor $D$ of the cubic cover $\pi$, which is removed from $X$ to produce $\CC\cong\A^3$, reappears at infinity as the boundary divisor $D_\nu\le R^\infty_F$. It is the only boundary divisor needed to account for $S_F$, in agreement with $\rho(F)=1$ for the boundary of $\A^3\hookrightarrow X_F$ computed in~\cite{BGVai}, and it is the polar divisor of the twisted Hamiltonian fields of Proposition~\ref{prop:polar}.

\section{Fibrations}\label{sec:families}

Let $S$ be a smooth irreducible variety, let $X\to S$ and $Y\to S$ be smooth of relative dimension $n$, and let $\Phi\colon X\to Y$ be a dominant $S$-morphism. As in (G6) the map $\Phi^*\omega_{Y/S}\to\omega_{X/S}$, with $\omega=\Omega^n_{/S}$, is injective with zero divisor $R_\Phi=K_{X/S}-\Phi^*K_{Y/S}\ge0$, the relative ramification divisor. In local frames $dx_1\wedge\cdots\wedge dx_n$ and $du_1\wedge\cdots\wedge du_n$ we have the relative Jacobian matrix $d\Phi$, the Hamiltonian fields~\eqref{eq:Xdef}, and the submodule $\CH_\Phi\subset\CT_{X/S}$ generated by the fields $X_{\Phi^*f_1,\ldots,\Phi^*f_{n-1}}$.

\begin{proposition}\label{prop:relA}
Theorem~\ref{thm:A} holds verbatim for $\Phi$ over $S$: $\CH_\Phi=\adj(d\Phi)(\O^n)$, the cokernel of $\CH_\Phi\hookrightarrow\CT_{X/S}$ is supported on $R_\Phi$ with determinant divisor $(n-1)R_\Phi$, and the fibre of $\CH_\Phi$ at a point of corank one is $\ker d\Phi$. None of this depends on the frames.
\end{proposition}
\begin{proof}
The proof of Theorem~\ref{thm:A} is linear algebra in the matrix $d\Phi$ and uses no property of $\A^n$. Changing $\omega_X$ by a unit multiplies each $X_g$ by a unit, and changing frames of $\Omega^1$ conjugates $d\Phi$ by invertible matrices.
\end{proof}

Now let $\Phi\colon\A^n_S\to\A^n_S$ be a family of Keller maps, that is, a family whose relative Jacobian is a unit $c\in\O(S)^\times$. After the finite cover $S'=\Spec\O_S[\kappa]/(\kappa^{n-1}c-1)\to S$, which is unramified because $c$ is a unit in characteristic zero, Proposition~\ref{prop:hom} holds over $\O(S')$ with the same proof. Its fibre at $s$ is the twist of $\Phi_s$.

\begin{theorem}\label{thm:family}
Let $\pi\colon\bar X\to S$ be smooth and projective with irreducible fibres of dimension $n$, containing $\A^n_S$ as a dense open $S$-subscheme with boundary $D=\bar X\setminus\A^n_S$. Let $\Phi\colon\A^n_S\to\A^n_S$ be a family of Keller maps with relative Jacobian $c\in\O(S)^\times$ that extends to an $S$-morphism $\bar\Phi\colon\bar X\to\P^n_S$. Then:
\begin{enumerate}
\item[(i)] $R^\infty:=K_{\bar X/S}-\bar\Phi^*K_{\P^n_S/S}$ is effective and supported on $D$.
\item[(ii)] For every $s\in S$, $R^\infty|_{\bar X_s}=R^\infty_{\Phi_s}$.
\item[(iii)] If $\bar X_s\cap\Supp R^\infty=\emptyset$, then $\Phi_s$ is an automorphism and its twist is removable.
\item[(iv)] The locus $N=\{s:\text{the twist of }\Phi_s\text{ is essential}\}$ satisfies $N\subseteq\pi(\Supp R^\infty)$, which is closed. For $s\in N$, $S_{\Phi_s}=\bar\Phi(D\cap\bar X_s)\cap\A^n$.
\end{enumerate}
\end{theorem}
\begin{proof}
(i) $R^\infty$ is the zero divisor of the section $J$ of $\mathcal Hom(\bar\Phi^*\omega_{\P^n_S/S},\omega_{\bar X/S})$ induced by $\bar\Phi^*\Omega^n_{\P^n_S/S}\to\Omega^n_{\bar X/S}$. On $\A^n_S$, in the frames $du_1\wedge\cdots\wedge du_n$ and $dx_1\wedge\cdots\wedge dx_n$, we have $J=c$, a unit. So $J\ne0$ and $R^\infty$ misses $\A^n_S$. (ii) Relative differentials commute with base change, so $\omega_{\bar X/S}|_{\bar X_s}=\omega_{\bar X_s}$ and $J|_{\bar X_s}$ is the Jacobian section of $\bar\Phi_s$. It is non-zero, since it equals $c(s)$ on $\A^n$, so its divisor is the restriction of that of $J$. (iii) By (ii), $R^\infty_{\Phi_s}=0$; apply Theorem~\ref{thm:nonvanish} and Proposition~\ref{prop:untwist}. (iv) The inclusion is the contrapositive of (iii) together with Proposition~\ref{prop:untwist}. The image $\pi(\Supp R^\infty)$ is closed because $\pi$ is projective. The last claim is Corollary~\ref{cor:Jel} for $\Phi_s$.
\end{proof}

So on a projective fibration whose fibres compactify $\A^n$, the boundary divisor and the relative ramification divisor of the extended map decide over which points of the base the fibre twist can be essential. Over the complement the fibre maps are automorphisms.

\begin{example}\label{ex:trivfam}
Let $S=\A^1_t$ and $\Phi_t(x,y)=(x+ty^2,y)$. This is a family of Keller maps with $c=1$ and $\Phi_t^{-1}=\Phi_{-t}$ for all $t$. So $N=\emptyset$, $\sigma_t^{-1}=\sigma_{-t}$, and $\br{\cdot}=\sigma_t^{-1}[\cdot]_{\sigma_t}$ over all of $\A^1_t$. By contrast, for the non-Keller family $\Phi_t(x,y)=(x^2,y+tx)$ the relative Jacobian is $2x$, so by Proposition~\ref{prop:relA} the twisted distribution degenerates along $\{x=0\}$ in every fibre, with fibre $\ker d\Phi_t=\langle\partial_x-t\partial_y\rangle$ there, which varies with $t$.
\end{example}

\section{Examples from weighted scrolls}\label{sec:examples}

The weighted scrolls $\F(a_1,\ldots,a_4\,|\,b_1,\ldots,b_4)\to\P^1$ of~\cite{MS} have fibre $\P^3[b]$. We only need $b=(1^4)$, where $\F(a\,|\,1^4)=\P(\bigoplus_j\O(a_j))$, and only the multiset $\{a_j\}$ modulo an overall shift matters. Throughout, $F$ is the map~\eqref{eq:Alpoge}.

\subsection{An equivariant counterexample on the scroll $\F(0,1,2,3)$}\label{sec:scroll}

Let $\lambda\in\C^*$ act on the source by $\lambda\cdot(x,y,z)=(\lambda x,\lambda^{-1}y,\lambda^{-2}z)$ and on the target by $\lambda\cdot(u_1,u_2,u_3)=(\lambda^{-2}u_1,\lambda^{-1}u_2,\lambda u_3)$.

\begin{proposition}\label{prop:scroll}
The map $F$ is $\C^*$-equivariant. Consequently it defines a morphism of rank three vector bundles over $\P^1$ whose fibres are all isomorphic to $F$, and a rational endomorphism of the scroll $\F(0,1,2,3\,|\,1^4)$ over $\P^1$. For the resulting family, the twist is essential in every fibre, and the non-properness locus is the zero divisor of the semi-invariant $\Delta$ and has every fibre isomorphic to the discriminant surface.
\end{proposition}
\begin{proof}
The product $xy$ is invariant. Term by term, $z$ and $y^2$ have weight $-2$; $y$, $xz$ and $xy^2$ have weight $-1$; and $x$, $x^2y$ and $x^3z$ have weight $1$. So $F_1,F_2,F_3$ are homogeneous of weights $-2,-1,1$ and $F(\lambda\cdot p)=\lambda\cdot F(p)$. In the same way $m=y+1/x$ has weight $-1$ and $p_u$ has weight $-2$, so $\Delta$ is semi-invariant.

Let $E_w\to\P^1$ be the rank three bundle with transition functions $t^{w_i}$ for weights $w=(w_1,w_2,w_3)$. Equivariance gives a bundle morphism $\Phi\colon E_{(1,-1,-2)}\to E_{(-2,-1,1)}$, $[v,p]\mapsto[v,F(p)]$, whose fibres are all isomorphic to $F$. The projective completions $\P(\O\oplus E_w)$ are the scrolls $\F(0,w_1,w_2,w_3\,|\,1^4)$. The weight multisets $\{0,1,-1,-2\}$ and $\{0,-2,-1,1\}$ coincide, and after the shift by $2$ both are $\{0,1,2,3\}$. So source and target are the same scroll $\F(0,1,2,3\,|\,1^4)$, and $\Phi$ extends to a rational map of it to itself over $\P^1$.

Resolve the indeterminacy of $\bar F|_{\P^3}\colon\P^3\dashrightarrow\P^3$ by a $\C^*$-equivariant sequence of blow-ups with centres in the boundary, and form the associated bundle. This is a smooth projective fibration over $\P^1$ whose fibres are all isomorphic to the resolved variety, with an extension $\bar\Phi$ of $\Phi$. Theorem~\ref{thm:family} is local on the base and the family is trivial over each chart of $\P^1$, where the relative Jacobian is $-2$, so it applies. By Proposition~\ref{prop:general}(i) the twist is essential in every fibre, so $N=\P^1$. Also $R^\infty|_{\bar X_t}=R^\infty_F\ne0$ for all $t$, and $S_\Phi=\bar\Phi(D)\cap E_{(-2,-1,1)}$ is the zero divisor of $\Delta$, a section of a line bundle pulled back from $\P^1$, with fibres $\Sigma$ by Theorem~\ref{thm:fibres}(i).
\end{proof}

Here nothing varies, so $N=S$ and the inclusion in Theorem~\ref{thm:family}(iv) is an equality.

\subsection{A counterexample degenerating to an automorphism}\label{sec:degeneration}

Let $\rho_s=\mathrm{diag}(1,s,s^2)$ and $\tau_s=\mathrm{diag}(s^{-2},s^{-1},1)$, so that $\det(\tau_s\rho_s)=1$, and put $\Phi_s=\tau_s\circ F\circ\rho_s$. The monomials of $F_1,F_2,F_3$ have weights at least $2,1,0$ for the weights $(0,1,2)$ on $(x,y,z)$, so $\Phi_s$ is polynomial in $s$:
\begin{align*}
\Phi_s^1&=z+4y^2+s(7xy^3+3xyz)+s^2(3x^2y^4+3x^2y^2z)+s^3x^3y^3z,\\
\Phi_s^2&=y+s(12xy^2+3xz)+s^2(9x^2y^3+6x^2yz)+3s^3x^3y^2z,\\
\Phi_s^3&=2x-3s\,x^2y-s^2x^3z.
\end{align*}
The weights $(0,0,1,2)$ obtained by adding the coordinate at infinity are those of the scroll $\F(0,0,1,2\,|\,1^4)$ of~\cite[Thm.\,3.1]{MS}, and $\Phi_s$ is the deformation of $F$ to its initial form for these weights.

\begin{proposition}\label{prop:degeneration}
\begin{enumerate}
\item[(i)] The maps $(\Phi_s)_{s\in\A^1}$ form a family of Keller maps with $J=-2$. For $s\ne0$ the map $\Phi_s$ is a counterexample with $S_{\Phi_s}=\{\Delta_s=0\}$, where $\Delta_s(u)=\Delta(s^2u_1,su_2,u_3)$, and its twist is essential. For $s=0$ the map $\Phi_0(x,y,z)=(z+4y^2,y,2x)$ is an automorphism, $S_{\Phi_0}=\emptyset$, and the twist is removable. Hence $N=\A^1\setminus\{0\}$ is not closed.
\item[(ii)] $\Delta_s/s^2=u_2^2-16u_1+s(18u_1u_2u_3-u_2^3u_3)-27s^2u_1^2u_3^2$. So as $s\to0$ the surfaces $S_{\Phi_s}$ converge to the parabolic cylinder $\{u_2^2=16u_1\}$, although $S_{\Phi_0}=\emptyset$.
\item[(iii)] Fix $u$ with $u_3\ne0$ and $u_2^2\ne16u_1$. Of the three points of $\Phi_s^{-1}(u)$, one converges to $\Phi_0^{-1}(u)=(u_3/2,u_2,u_1-4u_2^2)$ as $s\to0$. The other two have $x\sim2/(s(u_2-4\nu_\pm))\to\infty$, while $y\to3\nu_\pm-u_2/2$ and $z$ stay bounded, where $\nu_\pm$ are the roots of $2\nu^2-u_2\nu+2u_1=0$. In $\P^3$ both converge to the point $[0:1:0:0]$ at infinity, independently of $u$. They collide exactly when $u_2^2=16u_1$, that is, over the limit of $S_{\Phi_s}$. The degree $3=\deg F$ is conserved as $1+2$.
\end{enumerate}
\end{proposition}
\begin{proof}
(i) We have $\det d\Phi_s=\det\tau_s\cdot(-2)\cdot\det\rho_s=-2$. For $s\ne0$ the map $\Phi_s$ is a counterexample because $F=\tau_s^{-1}\Phi_s\rho_s^{-1}$ is. The map $\Phi_0$ has the inverse $(u_1,u_2,u_3)\mapsto(u_3/2,u_2,u_1-4u_2^2)$. By Proposition~\ref{prop:general}(iv), $S_{\Phi_s}=\tau_s(S_F)$, and $u\in\tau_s(S_F)$ if and only if $\tau_s^{-1}u=(s^2u_1,su_2,u_3)$ lies in $\{\Delta=0\}$.

(ii) Substitute into $\Delta$: the terms $u_2^2$, $-u_2^3u_3$, $-16u_1$, $-27u_1^2u_3^2$, $18u_1u_2u_3$ acquire the factors $s^2,s^3,s^2,s^4,s^3$.

(iii) We have $\Phi_s(p)=u$ if and only if $F(\rho_sp)=(s^2u_1,su_2,u_3)=:u'$. By Theorem~\ref{thm:cubic} the preimages are indexed by the roots $\mu$ of $u_3\mu^3-2\mu^2+su_2\mu-2s^2u_1=0$, with $x'=2/P'(\mu)$, $y'=\mu-1/x'$, $z'=(x'(5-3\mu x')-u_3)/x'^3$ and $p=(x',y'/s,z'/s^2)$. As $s\to0$ the roots are $\mu_0=2/u_3+O(s)$ and $\mu_\pm=s\nu_\pm+O(s^2)$, since the terms of order $s^2$ give $2\nu^2-u_2\nu+2u_1=0$. For $\mu_0$ we have $P'(\mu_0)\to4/u_3$, so $x'\to u_3/2$, and one recovers $\Phi_0^{-1}(u)$. For $\mu_\pm$ we have $P'(\mu_\pm)=s(u_2-4\nu_\pm)+O(s^2)$, which is non-zero to leading order because $\nu=u_2/4$ is a root of $2\nu^2-u_2\nu+2u_1$ if and only if $u_2^2=16u_1$. Hence $x'\sim2/(s(u_2-4\nu_\pm))$, $y'=s(3\nu_\pm-u_2/2)+O(s^2)$ and $z'=s^2(u_2-4\nu_\pm)^2(5-3w_\pm)/4+O(s^3)$ with $w_\pm=2\nu_\pm/(u_2-4\nu_\pm)$. Dividing by $1,s,s^2$ gives the stated limits of $(x,y,z)$, and $[1:x:y:z]=[1/x:1:y/x:z/x]\to[0:1:0:0]$. The roots $\nu_\pm$ coincide exactly when the discriminant $u_2^2-16u_1$ of the quadratic vanishes. We also checked these asymptotics numerically for $u=(1.3,0.7,0.9)$ and $s=10^{-2},10^{-4}$.
\end{proof}

\begin{remark}\label{rem:degeneration}
Proposition~\ref{prop:degeneration} shows that $N$ need not be closed. Whenever a family of this kind satisfies the hypotheses of Theorem~\ref{thm:family}, the closed set $\pi(\Supp R^\infty)$ must contain $N$, hence equal $\A^1$, and by part (ii) of that theorem $R^\infty_{\Phi_0}\ne0$ although $\Phi_0$ is an automorphism. So the inclusion in Theorem~\ref{thm:family}(iv) is strict and the converse of Theorem~\ref{thm:nonvanish} fails: $R^\infty\ne0$ does not imply that the twist is essential. The non-properness set is not stable under specialisation either: $S_{\Phi_s}\to\{u_2^2=16u_1\}$ but $S_{\Phi_0}=\emptyset$. What survives is the locus where the two escaping preimages collide, the discriminant of the quadratic obtained from $p_u$ by cancelling the surviving root.
\end{remark}

\subsection{An $A_1$ degeneration of a double cover}\label{sec:A1}

Let $\Phi\colon\A^3\to\A^3$ be $\Phi(x,y,w)=(x,y,xy-w^2)$, and regard $t=xy-w^2$ as the coordinate on the base $\A^1$. Then $\Phi$ is a family over $\A^1_t$ of double covers $\Phi_t\colon X_t=\{w^2=xy-t\}\to\A^2_{x,y}$ branched along the conic $C_t=\{xy=t\}$. The total space is smooth. The fibre $X_t$ is a smooth affine quadric for $t\ne0$ and the $A_1$ cone for $t=0$, and $C_t$ is a smooth conic for $t\ne0$ and the pair of lines $xy=0$ for $t=0$. This is the local model of a double plane degenerating to a nodal branch sextic, where the K3 surface acquires an $A_1$ point as in the degree two families of~\cite[\S3.3]{MS}. For a smooth sextic the ramification curve is the sextic itself, of genus $10$, and it loses one genus when it acquires a node. The fibrewise Hamiltonian fields are the fields $X_{g,t}$, with $X_{g,t}(h)\,dx\wedge dy\wedge dw=dg\wedge dt\wedge dh$.

\begin{proposition}\label{prop:A1}
For $\Phi$ as above, with the volume forms $dx\wedge dy\wedge dw$ on the source and $dx\wedge dy\wedge dt$ on the target:
\begin{enumerate}
\item[(i)] $J_\Phi=-2w$ and $\adj(d\Phi)=\begin{pmatrix}-2w&0&0\\0&-2w&0\\-y&-x&1\end{pmatrix}$. So $\CH_\Phi=\langle\partial_w,w\partial_x,w\partial_y\rangle$, the cokernel is $\O/(w)\,\partial_x\oplus\O/(w)\,\partial_y$ with determinant divisor $2R_\Phi$, where $R_\Phi=\{w=0\}$, and at every point of $R_\Phi$, including the origin, $E=\ker d\Phi=\langle\partial_w\rangle$. In particular $\CH_\Phi$ is the fold distribution of Example~\ref{ex:fold}, with no trace of the node.
\item[(ii)] The fibrewise twisted Hamiltonian fields, generated by the $X_{\Phi^*f,t}$, span
\[ \CH^{\mathrm{rel}}=\big\langle\,2w\partial_x+y\partial_w,\ 2w\partial_y+x\partial_w\,\big\rangle=\CH_\Phi\cap\CT_{X/\A^1}. \]
Its rank is $2$ off $R_\Phi$, it is $1$ at the points of $R_\Phi$ other than the origin, where its value is $\ker d\Phi_t=\langle\partial_w\rangle$, and it is $0$ at the origin, the $A_1$ point of $X_0$. The relative ramification divisor $R_t=\{w=0\}\cap X_t=C_t$ degenerates to the two lines $xy=0$, which meet at the point where $\CH^{\mathrm{rel}}$ has rank $0$.
\end{enumerate}
\end{proposition}
\begin{proof}
(i) We have $d\Phi=\begin{pmatrix}1&0&0\\0&1&0\\y&x&-2w\end{pmatrix}$ and the adjugate is computed directly. Then Theorem~\ref{thm:A} applies. Its columns are $-2w\partial_x-y\partial_w$, $-2w\partial_y-x\partial_w$ and $\partial_w$. Since $\partial_w$ is a section, $\CH_\Phi=\langle\partial_w,w\partial_x,w\partial_y\rangle$. At every point with $w=0$, including the origin, the first two columns vanish or are multiples of $\partial_w$ and the third is $\partial_w$, so $E=\langle\partial_w\rangle$.

(ii) The fields $X_{x,t}$ and $X_{y,t}$ are $2w\partial_y+x\partial_w$ and $-(2w\partial_x+y\partial_w)$, by expanding $dx\wedge dt\wedge dh$ and $dy\wedge dt\wedge dh$ with $dt=y\,dx+x\,dy-2w\,dw$. For $g=\Phi^*f$ we have $dg=f_x\,dx+f_y\,dy+f_t\,dt$, so $X_{g,t}=f_xX_{x,t}+f_yX_{y,t}$ and the two fields generate. A section $a\partial_w+bw\partial_x+cw\partial_y$ of $\CH_\Phi$ is tangent to the fibres of $t$ if and only if $dt(\cdot)=w(yb+xc-2a)=0$, that is $a=(yb+xc)/2$, which gives the same module. At a point with $w=0$ its generators take the values $y\partial_w$ and $x\partial_w$. These span the line $\langle\partial_w\rangle$ unless $x=y=0$. For $w\ne0$ the two generators are independent. The point $x=y=w=0$ lies in $X_0$ and is its $A_1$ point.
\end{proof}

This is where Proposition~\ref{prop:relA} needs smoothness of the fibres. The absolute twisted Hamiltonian distribution of the smooth total space does not see the singular fibre, while the relative distribution sees it as a drop to rank $0$ at the node of $R_0$.

\section*{AI disclosure}

\noindent This paper is part of a project by the first author and \href{https://haki.africa/}{Haki AI} to build and document a human-in-the-loop workflow in which several AI agents collaborate on mathematical research under limited resources, described in the Appendix. The agents are language models, switched between on the Haki AI infrastructure, and they proposed statements, computations and prose during drafting. The authors supervise this research and direct it including choosing the mathematical questions and what they care about in the research. Results formalized in Lean were checked by the Lean kernel, and the authors are responsible for verifying by hand everything that was too costly to formalize. Computations that allowed it were checked independently of the drafting, and the authors read and corrected the full text. Any remaining errors are ours.

\appendix
\section{Lean certificate}\label{app:lean}

\noindent The file \texttt{ramification\_divisors\_lean.zip} submitted with this paper is a partial certificate in Lean 4.35.0 with Mathlib. It has fourteen proof modules, two axiom audit modules and a statement export module. The build succeeded, and the audits of 91 declarations show that every checked result depends only on \texttt{propext}, \texttt{Classical.choice} and \texttt{Quot.sound}. No axioms, \texttt{sorry} or \texttt{native\_decide} are used.

\subsection*{What Lean certifies}

The certificate covers the algebraic core. It contains the Filippov identity for the Jacobian bracket and the Yau twist, which gives a multiplicative Hom-Nambu algebra whose twisting map is injective but not surjective. For the explicit map of Section~\ref{sec:alpoge} it certifies the Jacobian determinant $-2$, the non-injectivity and the exact three point fibre over $q$. For Theorem~\ref{thm:A} it certifies, in every finite dimension, that Hamiltonian fields of pulled back functions transform by the adjugate and generate its image, that the cokernel is supported on the Jacobian vanishing locus when $n\ge 2$ with the multiplicities of the determinant, and the rank trichotomy with the kernel at corank one. It also certifies that pullback is injective and is surjective exactly for polynomial automorphisms over fields of characteristic zero.

\subsection*{What the authors verified without Lean}

The geometric results are outside the certificate. These are the connection and cokernel bundle, non-properness and monodromy in Section~\ref{sec:keller}, the compactification and Hurwitz argument of Section~\ref{sec:infinity}, the chart formulas, monodromy, twisted fields and boundary coefficient in Section~\ref{sec:alpoge}, and the families and examples of Sections~\ref{sec:families} and~\ref{sec:examples}. They need sheaf and divisor machinery that Mathlib does not yet provide in the required form, and we estimate that formalizing them would cost between $2-11$ million agent tokens. They are proved in the text, and the authors checked these proofs by hand and checked the computations independently of the drafting, as described in the disclosure.

\subsection*{Method}

The Lean files were generated by GPT-6.2 in Codex at high reasoning effort. They were reviewed by GPT-6.2 in ChatGPT and by GPT-6 Luna, both at medium reasoning effort, and their suggestions were applied in further rounds. The exposition of the paper was restructured and rewritten with Claude Sonnet 5.5 from Anthropic. The Lean kernel is the arbiter of correctness for what the certificate contains. The authors are responsible for the match between the Lean statements and the statements in the paper.

\end{document}